
\input amstex.tex
\documentstyle{amsppt}
\magnification1200

\hsize=12.5cm \vsize=18cm \hoffset=1cm \voffset=2cm

\footline={\hss{\vbox to 2cm{\vfil\hbox{\rm\folio}}}\hss}
\nopagenumbers
\def\DJ{\leavevmode\setbox0=\hbox{D}\kern0pt\rlap
{\kern.04em\raise.188\ht0\hbox{-}}D}

\def\txt#1{{\textstyle{#1}}}
\baselineskip=13pt
\def\hf{{\textstyle{1\over2}}}
\def\a{\alpha}\def\b{\beta}
\def\d{{\,\roman d}}
\def\e{\varepsilon}
\def\f{\varphi}
\def\G{\Gamma}
\def\k{\kappa}
\def\s{\sigma}

\def\={\;=\;}

\def\zt{\zeta(\hf+it)}

\def\D{\Delta}

\def\R{\Re{\roman e}\,} \def\I{\Im{\roman m}\,}
\def\z{\zeta}

\def\H{H_j^3({\txt{1\over2}})} 
\def\hf{{\textstyle{1\over2}}}
\def\txt#1{{\textstyle{#1}}}
\def\f{\varphi}

\font\tenmsb=msbm10
\font\sevenmsb=msbm7
\font\fivemsb=msbm5
\newfam\msbfam
\textfont\msbfam=\tenmsb
\scriptfont\msbfam=\sevenmsb
\scriptscriptfont\msbfam=\fivemsb
\def\Bbb#1{{\fam\msbfam #1}}

\def \NN {\Bbb N}

\def \RR {\Bbb R}
\def \ZZ {\Bbb Z}

\font\ff=cmr8
\def\txt#1{{\textstyle{#1}}}
\baselineskip=13pt

\font\teneufm=eufm10
\font\seveneufm=eufm7
\font\fiveeufm=eufm5
\newfam\eufmfam
\textfont\eufmfam=\teneufm
\scriptfont\eufmfam=\seveneufm
\scriptscriptfont\eufmfam=\fiveeufm
\def\mathfrak#1{{\fam\eufmfam\relax#1}}

\font\tenmsb=msbm10
\font\sevenmsb=msbm7
\font\fivemsb=msbm5
\newfam\msbfam
     \textfont\msbfam=\tenmsb
      \scriptfont\msbfam=\sevenmsb
      \scriptscriptfont\msbfam=\fivemsb

  \def\rightheadline{{\hfil{\ff
  Moments of $|\zt|$ in short intervals}\hfil\tenrm\folio}}

  \def\leftheadline{{\tenrm\folio\hfil{\ff
   A. Ivi\'c }\hfil}}
  \def\emptyheadline{\hfil}
  \headline{\ifnum\pageno=1 \emptyheadline\else
  \ifodd\pageno \rightheadline \else \leftheadline\fi\fi}

\topmatter
\title
ON MOMENTS OF $|\zt|$ IN SHORT INTERVALS
\endtitle
\author   Aleksandar Ivi\'c  \endauthor
\address
Aleksandar Ivi\'c, Katedra Matematike RGF-a
Universiteta u Beogradu, \DJ u\v sina 7, 11000 Beograd,
Serbia (Yugoslavia).
\endaddress
\keywords
 Riemann zeta-function, mean
square and fourth moment of $|\zt|$, moments of $|\zt|$ in short intervals
\endkeywords
\subjclass
11M06 \endsubjclass
\email {\tt
aivic\@rgf.bg.ac.yu,  aivic\@matf.bg.ac.yu} \endemail
\dedicatory
Dedicated to Prof. K. Ramachandra on the occasion of his seventieth
birthday
\enddedicatory
\abstract
{Power moments of
$$
J_k(t,G) = {1\over\sqrt{\pi}G}
\int_{-\infty}^\infty |\z(\hf + it + iu)|^{2k}{\roman e}^{-(u/G)^2}\d u
\qquad(t \asymp T, T^\e \le G \ll T),$$
where $k$ is a natural number,
are investigated. The results that are obtained are used to show
how bounds for $\int_0^T|\zt|^{2k}\d t$ may be obtained.
}
\endabstract
\endtopmatter

\head
1. Introduction
\endhead
Power moments represent one of the most important parts of the
theory of the Riemann zeta-function $\z(s) = \sum_{n=1}^\infty n^{-s}$
$\;(\s = \R s > 1)$. Of particular significance are the moments on the
``critical line" $\s = \hf$, and a vast literature exists on this
subject (see e.g., [8], [9], [20], [22], [24] and [26]). Let us define
$$
I_k(T) = \int_0^T|\zt|^{2k}\d t,\leqno(1.1)
$$
where $k\in\RR$ is a fixed, positive number. Naturally one would want to
find an asymptotic formula for $I_k(T)$ for a given $k$, but this is an
extremely difficult problem. Except when $k = 1$ and $k = 2$, no asymptotic
formula for $I_k(T)$ is known yet, although there are plausible conjectures
for such formulas (see e.g., [2]). In the absence of asymptotic
formulas for $I_k(T)$, one would like then to obtain upper and lower
bounds for $I_k(T)$, and for the closely related problem of
$$
I_k(T+G) - I_k(T-G) = \int_{T-G}^{T+G}|\zt|^{2k}\d t\quad(1\ll G \le T).
\leqno(1.2)
$$
For the latter, important results were obtained by K. Ramachandra,
either alone, or in collaboration with R. Balasubramanian. Many of his
results are contained in his comprehensive monograph [24] on mean
values and omega-results for the Riemann zeta-function. In particular,
[24] contains the proof of the lower bound
$$
\int_{T-G}^{T+G}|\zt|^{2k}\d t \gg_k G(\log G)^{k^2}
\quad(\log\log T \ll_k G \le T, k\in \NN),\leqno(1.3)
$$
where $\ll_k$ (or $\gg_k$) means that the implied constant depends only
on $k$. One believes that the bound in (1.3) represents the correct order
 of magnitude, at least for a certain range of $G$ for a given $k\in \NN$.
Unfortunately, even proving the corresponding much weaker upper bound
(for $G=T$), namely
$$
I_k(T) \;\ll_{\e,k} T^{1+\e}\qquad(k > 0)\leqno(1.4)
$$
seems at present impossible for any $k>2$. Here and later, $\e >0$ denotes
constants which may be arbitrarily small, but are not necessarily the
same ones at each occurrence. In view of the relation
(see  [9] or [24])
$$
\z^k(\hf+it) \ll_k \log t\left(\int_{t-1/3}^{t+1/3}|\z(\hf+iu)|^k \d t
\right) + 1\qquad(k\in\NN),\leqno(1.5)
$$
it is easily seen that (1.4) (for all $k$) is equivalent to the famous
{\it Lindel\"of hypothesis} that $\zt \ll_\e |t|^\e$. The Lindel\"of
hypothesis, like the even more famous {\it Riemann hypothesis} (that
all complex zeros of $\z(s)$ have real part 1/2), is neither proved
nor disproved at the time of writing of this text. For a discussion
on this subject, see [13].

\medskip
The aim of this paper is to investigate upper bounds for $I_k(T)$
when $k\in\NN$, which we henceforth assume. The problem can be reduced
to bounds of $|\zt|$ over short intervals, as in (1.2), but it
is more expedient to work with the smoothed integral
$$
J_k(T,G) := {1\over \sqrt{\pi}G}\int_{-\infty}^\infty
|\z(\hf + iT + iu)|^{2k}{\roman e}^{-(u/G)^2}\d u\quad(1 \ll G \ll T).
\leqno(1.6)
$$
Namely we obviously have
$$
I_k(T+G) - I_k(T-G) = \int_{-G}^G|\z(\hf + iT + iu)|^{2k}\d u
\le \sqrt{\pi}{\roman e}G\,J_k(T,G),\leqno(1.7)
$$
and it is technically more convenient to work with $J_k(T,G)$ than with
$I_k(T+G) - I_k(T-G)$. Of course, instead of the Gaussian exponential
weight $\exp(-(u/G)^2)$, one could introduce in (1.6) other smooth weights
with a similar effect. The Gaussian weight has the advantage that, by the
use of the classical integral
$$
\int_{-\infty}^\infty \exp(Ax-Bx^2)\d x = \sqrt{\pi\over B}\,\exp\left(
{A^2\over4B}\right)\qquad(\R B > 0),\leqno(1.8)
$$
one can often explicitly evaluate the relevant exponential integrals
that appear in the course of the proof.

\medskip
The plan of the paper is as follows. In the next two sections we
shall briefly discuss the results on $I_k(T)$ and $J_k(T,G)$ when
$k= 1$ and $k=2$, respectively. Indeed, as these are the only cases
when we possess relatively good knowledge and explicit formulas, it
is only natural that those results be used in deriving results on
higher power moments, when our knowledge is quite imperfect. We shall obtain
new results on moments of $J_k(T,G)$ by using the explicit
formulas of Section 2 and Section 3. This will be
done in Section 4 and Section 5. Finally, in Section 6, it will be
shown how one can obtain bounds for $I_k(T)$ from the bounds of moments
of $J_k(T,G)$.

\head
2. The mean square formula
\endhead
The mean square formula for $|\zt|$ is traditionally written in the form
$$
\int_0^T|\zt|^2\d t \= T\log\left({T\over2\pi}\right) + (2\gamma
- 1)T + E(T),\leqno(2.1)
$$
where $\gamma = -\G'(1) = 0.577\ldots\,$ is Euler's constant, and
$E(T)$ is to be considered as the  error term in the
asymptotic formula (2.1). F.V. Atkinson [1] established in 1949 an explicit,
albeit complicated formula for $E(T)$, containing two exponential sums of
length $\asymp T$ weighted by the number of divisors function $d(n)$,
plus an error term which is $O(\log^2T)$. This is given as

\medskip
LEMMA 1. {\it Let $0 < A < A'$ be any two fixed constants
such that $AT < N < A'T$, and let $N' = N'(T) =
T/(2\pi) + N/2 - (N^2/4+ NT/(2\pi))^{1/2}$. Then }
$$
E(T) = \Sigma_1(T) + \Sigma_2(T) + O(\log^2T),\leqno(2.2)
$$
{\it where}
$$
\Sigma_1(T) = 2^{1/2}(T/(2\pi))^{1/4}\sum_{n\le N}(-1)^nd(n)n^{-3/4}
e(T,n)\cos(f(T,n)),\leqno(2.3)
$$
$$
\Sigma_2(T) = -2\sum_{n\le N'}d(n)n^{-1/2}(\log T/(2\pi n))^{-1}
\cos(T\log T/(2\pi n) - T + \pi /4),\leqno(2.4)
$$
{\it with}
$$
\eqalign{\cr&
f(T,n) = 2T{\roman {arsinh}}\,\bigl(\sqrt{\pi n/(2T})\bigr) + \sqrt{2\pi nT
+ \pi^2n^2} - \txt{1\over4}\pi \cr&
=  -\txt{1\over4}\pi + 2\sqrt{2\pi nT} +
\txt{1\over6}\sqrt{2\pi^3}n^{3/2}T^{-1/2} + a_5n^{5/2}T^{-3/2} +
a_7n^{7/2}T^{-5/2} + \ldots\,,\cr}\leqno(2.5)
$$
$$\eqalign{\cr
e(T,n) &= (1+\pi n/(2T))^{-1/4}{\Bigl\{(2T/\pi n)^{1/2}
{\roman {arsinh}}\,\bigl(\sqrt{\pi n/(2T)}\,\bigr)\Bigr\}}^{-1}\cr&
= 1 + O(n/T)\qquad(1 \le n < T),
\cr}\leqno(2.6)
$$
{\it and $\,{\roman{arsinh}}\,x = \log(x + \sqrt{1+x^2}\,).$}

Atkinson's formula was the starting point for many results
on $E(T)$ (see [8, Chapter 15] for some of them).
It is conjectured that $E(T) \ll_\e T^{1/4+\e}$, but currently
this bound cannot be proved even if the Riemann Hypotheis is assumed. The
best known upper bound for $E(T)$, obtained by intricate estimation of
a certain exponential sum, is due to M.N. Huxley [6]. This is
$$
E(T) \;\ll\; T^{72/227}(\log T)^{679/227}, \quad {72\over227} =
0.3171806\ldots\;.
$$
In the other direction, J.L. Hafner and the author [3] proved
 that there exist absolute constants $A, B > 0$ such  that
$$
E(T) = \Omega_+\left\{(T\log T)^{1/4}(\log\log T)^{(3+\log4)/4}
{\roman e}^{-A\sqrt{\log\log\log T}}\right\}
$$
and
$$
E(T) = \Omega_-\left\{T^{1/4}\exp\left({B(\log\log T)^{1/4}
\over(\log\log\log T)^{3/4}}\right)\right\},
$$
where $f(x) = \Omega_+(g(x))$ means that
$\limsup\limits_{x\to\infty} f(x)/g(x) >  0,$  and
$f(x) = \Omega_-(g(x))$ means that
$\liminf\limits_{x\to\infty} f(x)/g(x) < 0$.

In what follows we shall formulate an explicit formula for
$J_1(T,G)$. Such a result can be, of course, deduced from
Atkinson's formula (2.3)-(2.5) by the use of (1.8). This approach was used
originally by D.R. Heath-Brown  [4], who proved
$$
\int_0^T|\zt|^{12}\d t \ll T^2\log^{17}T,\leqno(2.7)
$$
which is still essentially the best known result concerning higher
power moments of $|\zt|$. This procedure can be avoided by appealing to
Y. Motohashi's formula [22, p. 213], which states that
$$
J_1(T,G) = 2^{3\over4}\pi^{1\over4}T^{-{1\over4}}\sum_{n=1}^\infty
(-1)^nd(n)n^{-{1\over4}}\sin f(T,n)
\exp\Bigl(-{\pi nG^2\over2T}\Bigr) + O(\log T),
\leqno(2.8)
$$
where $f(T,n)$ is given by (2.5), and $T^{1/4} \le G \le T/\log T$.
In fact, only the range $G \le T^{1/3}$ is relevant, since for
$G \ge T^{1/3}$ one has $J_1(T,G) \ll \log T$ by [8, Chapter 7].
Motohashi's proof of (2.8), like the proof of Atkinson's formula
for $E(T)$, is based on classical methods from analytic number theory.
Albeit the expression on the right-hand side of (2.8) is quite simple,
the condition $G \ge T^{1/4}$ is rather restrictive for the application that
we have in mind. Thus we shall use a similar type of result, which is
valid in a much wider range. This is contained in

\medskip
LEMMA 2. {\it For $T^\e \le G \le T$ and $f(T,n)$ given by} (2.5),
{\it we have}
$$
\eqalign{&
J_1(T,G) = O(\log T) \;+ \cr&
+ \sqrt{2}\sum_{n=1}^\infty (-1)^nd(n)n^{-1/2}
\left({\Bigl({T\over2\pi n} + {1\over4}\Bigr)}^{1/2}
- {1\over2}\right)^{-1/2}\times\cr&\times\exp\left(-G^2(
{\roman {arsinh}}\,\sqrt{\pi n/(2T)})^{2}\right)\sin f(T,n).\cr}
\leqno(2.9)
$$

\bigskip
By using Taylor's formula it is seen that the
error made by replacing
$$
\left(\left({T\over2\pi n} + {1\over4}\right)^{1/2}
- {1\over2}\right)^{-1/2}
\exp\left(-G^2({\roman {arsinh}}\,\sqrt{\pi n/(2T)})^{2}\right)
$$
by
$$
\left({T\over2\pi n}\right)^{-1/4}\exp(-\pi G^2/(nT))
$$
is $\ll 1$ for
$G \ge T^{1/5}\log^CT$. But the important fact is that in applications
(2.9) is as useful as (2.8), since the factors under the sine function
are identical.

\medskip
{\bf Proof of Lemma 2.} The proof of (2.9) follows fom Y. Motohashi
[22, Theorem 4.1], which gives that
$$
\eqalign{&
\int_{-\infty}^\infty |\zt|^2g(t)\d t = \int_{-\infty}^\infty
\left[\R\left\{{\G'\over\G}(\hf+it)\right\} + 2\gamma - \log(2\pi)\right]
g(t)\d t \cr&+ 2\pi\R(g(\hf i))
+ 4\sum_{n=1}^\infty d(n)\int_0^\infty (y(y+1))^{-1/2}g_c(\log(1+1/y))
\cos(2\pi ny)\d y,\cr}\leqno(2.10)
$$
where
$$
g_c(x) := \int_{-\infty}^\infty g(t)\cos(xt)\d t
$$
is the cosine Fourier transform of $g(t)$. One requires the function $g(r)$
to be real-valued for $r\in\RR$, and that there exists a large constant
$A>0$ such that $g(r)$ is regular and $\ll (|r|+1)^{-A}$ for
$|\I r| \le A$. The choice
$$
g(t) = {1\over\sqrt{\pi} G}{\roman e}^{-(T-t)^2/G^2},
\quad g_c(x) = {\roman e}^{-{1\over4}(Gx)^2}\cos(Tx)
$$
is permissible, and then the integral on the left-hand side of (2.10)
becomes $J_1(T,G)$. The first integral on the right-hand side of
(2.10) is $O(\log T)$, and the second one is evaluated by the
saddle-point method (see e.g., [8, Chapter 2]). A convenient result to
use is [8, Theorem 2.2 and Lemma 15.1], due originally to Atkinson [1].
 In the latter only the exponential
factor $\exp(-{1\over4}G^2\log(1+1/y))$ is missing. In the notation
of [1] and [8] we have that the saddle point $x_0$ satisfies
$$
x_0 = U - {1\over2} = \left({T\over2\pi n} +
{1\over4}\right)^{1/2} - {1\over2},
$$
and the presence of the above exponential factor makes it possible to
truncate the series in (2.9) at $n = TG^{-2}\log T$ with a negligible error.
Furthermore, in the remaining range for $n$ we have
$$
\Phi_0\mu_0F_0^{-3/2} \ll (nT)^{-3/4},
$$
which makes a total copntribution of $O(1)$, as does error term
integral in Theorem 2.2 of [8]. The error terms with $\Phi(a),\,\Phi(b)$
vanish for $a = 0,\,b=\infty\,$, and (2.9) follows.

\head
3. The  formula for the fourth moment
\endhead

The asymptotic formula for the fourth moment of the Riemann zeta-function
$\zeta(s)$ on the critical line is customarily written as
$$
\int_0^T|\zt|^4\d t \;=\; TP_4(\log T) \;+\;E_2(T),\;
P_4(x) \;=\; \sum_{j=0}^4\,a_jx^j.\leqno(3.1)
$$
A classical result of A.E. Ingham [7] from 1926 is that $a_4 = 1/(2\pi^2)$
and that the error term $E_2(T)$
in (3.1) satisfies the bound $E_2(T) \ll T\log^3T$ (a simple proof
of this is due to K. Ramachandra [23]). Much later D.R. Heath-Brown
[4] made  progress in this problem by proving that
$E_2(T) \ll_\e T^{7/8+\e}$. He also calculated
$$
a_3 \;=\; 2(4\gamma - 1 - \log(2\pi) - 12\zeta'(2)\pi^{-2})\pi^{-2}
$$
and produced more complicated expressions for $a_0, a_1$ and $a_2$ in (3.1).
For an explicit evaluation of the $a_j$'s  in (3.1) the reader
is referred to the author's work [11]. In the last fifteen years,
due primarily to the application of powerful methods of spectral theory
(see Y. Motohashi's monograph [22] for
a comprehensive account), much advance has been made in connection with
$E_2(T)$.  This involves primarily results with
exponential sums involving the quantities
$\kappa_j$ and $\a_jH_j^3(\hf)$. Here as usual
$\,\{\lambda_j = \k_j^2 + {1\over4}\} \,\cup\, \{0\}\,$
is the discrete spectrum of the non-Euclidean Laplacian acting
on $\,SL(2,\ZZ)\,$--automorphic forms, and
$\a_j = |\rho_j(1)|^2(\cosh\pi\k_j)^{-1}$, where
$\rho_j(1)$ is the first Fourier coefficient of the Maass wave form
corresponding to the eigenvalue $\lambda_j$ to which the Hecke series
$H_j(s)$ is attached. It is conjectured that $E_2(T) \ll_\e T^{1/2+\e}$,
which would imply the (hitherto unproved) bound $\zt \ll_\e |t|^{1/8+\e}$.
It is known now that
$$
E_2(T) \;=\; O(T^{2/3}\log^{C_1}T),\quad E_2(T) \;=\; \Omega(T^{1/2}),
\leqno(3.2) $$
$$
\int_0^TE_2(t)\d t \;=\; O(T^{3/2}),\quad\int_0^TE_2^2(t)\d t
\;=\; O(T^2\log^{C_2}T),\leqno(3.3)
$$
with effective constants $C_1,\,C_2 > 0$ (the values $C_1 = 8, C_2 = 22$ are
worked out in [22]). The above results were proved by
Y. Motohashi and the author: (3.2) and the first bound in (3.3) in
[17], and the second upper bound in (3.3) in [16]. The $\Omega$--result
in (3.2)
was improved to $E_2(T) = \Omega_\pm(T^{1/2})$ by Y. Motohashi [21].
It turns out that there is no explicit formula for $E_2(T)$ which would
represent the analogue of Atkinson's formula (cf. Lemma 1). Results
on $E_2(T)$ have been obtained indirectly, by using the explicit
formula for $J_2(T,G)$, due to Y. Motohashi (see [22]). This is

\bigskip
LEMMA 3. {\it Let $D>0$ be an arbitrary constant. For $T^{1/2}\log^{-D}T
\le G \le T/\log T$ we have}
$$\eqalign{
J_2(T,G) &= O(\log^{3D+9}T)  \cr&
+ {\pi\over\sqrt{2T}}\sum_{j=1}^\infty \a_j H_j^3(\hf)\k_j^{-1/2}
\sin\left(\k_j\log{\k_j\over4{\roman e}T}\right)
\exp(-{\txt{1\over4}}(G\k_j/T)^2).\cr}\leqno(3.4)
$$

\bigskip
In what concerns higher moments, let us only state that from (2.7)
and (3.1) one obtains by H\"olders's inequality for integrals
$$
\int_0^T|\zt|^6\d t \ll T^{5/4}\log^{29/4}T,
\quad \int_0^T|\zt|^8\d t \ll T^{3/2}\log^{21/2}T. \leqno(3.5)
$$
The bounds in (3.5) are hitherto the sharpest ones known.

Let it be also mentioned here that the sixth moment was investigated
by the author in [12], where it was shown that
$$
\int_0^T|\zt|^6\d t  \ll_\e T^{1+\e}\leqno(3.6)
$$
does hold if a certain conjecture involving the so-called ternary additive
divisor problem is true.

\head
4. The  moments of $J_1(t,G)$
\endhead
In this section we shall prove results on moments of $J_1(t,G)$. One
expects this function, at least for certain ranges of $G$, to behave
like $O(t^\e)$ on the average. Our bounds are contained in

\medskip
THEOREM 1. {\it We have}
$$
\int_T^{2T}J_1^m(t,G)\d t \ll_\e T^{1+\e}\leqno(4.1)
$$
{\it for $T^\e \le G \le T$ if $m = 1,2$; for $T^{1/7+\e} \le G \le T$ if
$m=3$, and for $T^{1/5+\e} \le G \le T$ if $m=4$.}

\medskip
{\bf Proof of Theorem 1.} Our starting point in all cases is the
explicit formula (2.9). The results for $m=1$ and $m=2$ follow by
straightforward integration and the first derivative test for
exponential integrals (see [8, Lemma 2.1]). The proof resembles
mean square bounds for $\D(x)$ (the error term in the divisor problem)
and $E(t)$ (op. cit.), and is omitted for the sake of brevity. Instead,
we shall concentrate on the more difficult cases $m=3$ and $m=4$.
For this we shall need two lemmas on the spacing of three and four
square roots (the square roots appear in view of the asymptotic
formula given in (2.5)). These are

\medskip
LEMMA 4. {\it Let ${\Cal N}$ denote the number of solutions  in integers
m,n,k of the inequality}
$$
|\sqrt{m} + \sqrt{n} - \sqrt{k}| \le \delta\sqrt{M} \qquad(\delta > 0)
$$
{\it with $M' < n \le 2M', M < m \le 2M$, and $M' \le M$. Then}
$$
{\Cal N} \ll_\e M^{\e}(M^2M'\delta + (MM')^{1/2}).\leqno(4.2)
$$

\medskip
LEMMA 5  {\it Let $k\ge 2$ be a fixed
integer and $\delta > 0$ be given.
Then the number of integers $n_1,n_2,n_3,n_4$ such that
$N < n_1,n_2,n_3,n_4 \le 2N$ and}
$$
|n_1^{1/k} + n_2^{1/k} - n_3^{1/k} - n_4^{1/k}| < \delta N^{1/k}
$$
{\it is, for any given $\e>0$,}
$$
\ll_\e N^\e(N^4\delta + N^2).\leqno(4.3)
$$

\medskip
Lemma 4 was proved by Sargos and the author [18], while Lemma 5 is due
to Robert--Sargos [25]. The plan of the proof of (4.1) when $m=3$
is simple: the expression (2.9) will be raised to the third power
and then integrated. There are, however, two obstacles in attaining
this goal. The first is that direct integration does not lead to
adequate truncation, so that some smoothing of the relevant
integral will be made. The second one is that in the asymptotic formula
for $f(T,n)$ in (2.5) not only square roots appear, but also higher
powers. To get around this difficulty we shall appeal to

\medskip
LEMMA 6 (M. Jutila [9]). {\it For $A\in\RR$ we have} $$
\cos\left(\sqrt{8\pi nT} +
\txt{1\over6}\sqrt{2\pi^3}n^{3/2}T^{-1/2} - A\right) =
\int_{-\infty}^\infty \a(u)\cos(\sqrt{8\pi n}(\sqrt{T} + u) - A)\d
u,\leqno(4.4) $$ {\it where $\a(u) \ll T^{1/6}$ for $u\not=0$,} $$
\a(u) \ll T^{1/6}\exp(-bT^{1/4}|u|^{3/2}) \leqno(4.5) $$ {\it for
$u<0$, and} $$ \a(u) =
T^{1/8}u^{-1/4}\left(d\exp(ibT^{1/4}u^{3/2}) + {\bar
d}\exp(-ibT^{1/4}u^{3/2})\right) + O(T^{-1/8}u^{-7/4})\leqno(4.6)
$$ {\it for $u \ge T^{-1/6}$ and some constants $b\; (>0)$ and
$d$.}

\medskip
Now we continue with the proof of Theorem 1. Write first $$
\int_T^{2T}J_1^m(t,G)\d t \le  \int_{T/2}^{5T/2}\f(t)J_1^m(t,G)\d
t, \leqno(4.7) $$ where $\f(t) \;(\ge 0)$ is a smooth function
supported in $[T/2,\,5T/2]\,$ such that $\f(t) = 1$ when $t\in
\,[T,\,2T]\,$, and then we have $\f^{(r)}(t) \ll_r T^{-r}\;(r =
0,1,2,\ldots)\,$. We truncate (2.9) at $TG^{-2}\log T$ and use it
to expand the $m$-th power on the right-hand side of (4.7) when $m
= 3,4$. The terms $a_{2j-2}n^{(2j-1)/2}T^{-j/2}\;(j\ge3)$  in
$f(T,n)$ (cf. (2.5)) are expanded by Taylor's formula. Since $$
n^{5/2}t^{-3/2} \ll TG^{-5}\log^{5/2}T \le T^{-\e}
 $$
 for $G \ge
T^{1/5+\e}$, it transpires that in this range for $G$ we may take
sufficiently many terms in Taylor's formula so that the error term
will make a negligible contribution. The other terms will lead to
similar expressions, and the largest contribution will come from
the constant term. After this there will remain $$ \sin h(t,k)\sin
h(t,\ell)\sin h(t,n),\; h(t,u) = \sqrt{8\pi tu} +
{\txt{1\over6}}\sqrt{2\pi^3}u^{3/2}t^{-1/2} - {\txt{1\over4}}\pi
$$
when $m=3$, with $T^{1/3} < k,\ell,n \le TG^{-2}\log T$, $k \asymp
K, k \ge \max(\ell,n)$. The factors with $u^{3/2}$ (for $u =
k,\ell,n$) are removed from the $h$-functions by the use of Lemma
6. With $\a(v)$ given by (4.6) we have $$ \eqalign{\cr&
\cos\left(\sqrt{8\pi nt} +
\txt{1\over6}\sqrt{2\pi^3}n^{3/2}t^{-1/2}- A\right) = O(T^{-10})\,
+ \cr&
 \int_{-u_0}^{u_1}\a(v)\cos(\sqrt{8\pi n}(\sqrt{t}+v) -
A)\d v  +
\int_{u_1}^\infty \a(v)\cos(\sqrt{8\pi n}(\sqrt{t}+v) -
A)\d v ,\cr}\leqno(4.8)
$$
where we set
$$
u_0 =  T^{-1/6}\log T, \; u_1 = CKT^{-1/2},\leqno(4.9)
$$
and $C>0$ is a large constant.  With this choice of $u_0,u_1$ and
(4.5)-(4.6) it follows that, for $T/2 \le t \le 5T/2$,
$$
 \int_{-u_0}^{u_1}\a(v)\cos(\sqrt{8\pi n}(\sqrt{t}+v) -
A)\d v  + \int_{u_1}^\infty \a(v)\cos(\sqrt{8\pi n}(\sqrt{t}+v) -
A)\d v \ll \log T.\leqno(4.10) $$ Namely we have $$
\int_{u_0}^{u_1}t^{1/8}v^{-1/4}\exp(ibt^{1/4}v^{3/2} \pm
\sqrt{8\pi n}v)\d v \ll \log T, $$ on writing the integral as a
sum of $\ll \log T$ integrals over $[U,\,U']$ with $u_0 \le U < U'
\le 2U \ll u_1$, and applying the second derivative test (i.e.,
[8, Lemma 2.2]) to each of these integrals. We remark that the
contribution of the $O$-term in (4.6) will be, by trivial
estimation, $O(1)$. It remains yet to deal with the integral with
$v > u_1$ in (4.10), when we note that $$ {\partial \over \partial
v} \left(bt^{1/4}v^{3/2} \pm \sqrt{8\pi n}v\right) \;\gg\;
T^{1/4}v^{1/2} \quad(v > u_1,\;t\asymp T), $$ provided that $C$ in
(4.9) is sufficiently large. Thus by the first derivative test $$
\eqalign{\cr& \int_{u_1}^\infty \a(v)\cos(\sqrt{8\pi
n}(\sqrt{t}+v) - A )\d v \cr& \ll 1 +
T^{1/8}u_1^{-1/4}T^{-1/4}u_1^{-1/2} \cr& \ll 1 + T^{1/4}K^{-3/4}
\ll 1, \cr} $$ since $K \gg T^{1/3}$. Thus (4.10) holds.

Hence setting $$ E_\pm \;:=\; \sqrt{8\pi}(\sqrt{k} + \sqrt{\ell}
\pm \sqrt{n}\,), $$ it is seen that we are left with the integral
of $$ \int_{T/2}^{5T/2}\f(t)F(t;k,\ell,n){\roman
e}^{iE_\pm\sqrt{t}}\d t, \leqno(4.11) $$ where $$ \eqalign{&
F(t;k,\ell,n) := \left(\left({t\over2\pi k}+{1\over4}\right)^{1/2}
- {1\over2}\right)^{-1/2} \left(\left({t\over2\pi
\ell}+{1\over4}\right)^{1/2} - {1\over2}\right)^{-1/2}\cr&
\times\left(\left({t\over2\pi n}+{1\over4}\right)^{1/2} -
{1\over2}\right)^{-1/2}\exp(-G^2({\roman {arsinh}}\sqrt{\pi
k/(2t)})^2) \cr& \times\exp(-G^2({\roman {arsinh}}\sqrt{\pi
\ell/(2t)})^2) \exp(-G^2({\roman {arsinh}}\sqrt{\pi
n/(2t)})^2).\cr} $$ Repeated integration by parts show that the
integral in (4.11) with $E_+$ will make a negligible contribution,
and also the one with
 $|E_-| \ge T^{\e-1/2}$ for any
given $\e >0$. The contribution of those $k,\ell,n$ for which
$|E_-| \le T^{\e-1/2}$ is estimated by the use of Lemma 4 (with an
obvious change of notation and with $\delta \asymp
K^{-1/2}T^{\e-1/2}$). After this, the integral over $t$ is
estimated trivially, and (4.10) is used. The relevant expression
on the right-hand side of (4.7) is $$ \eqalign{& \ll_\e
T^{1+\e}\max_{T^{1/3}\le K \le TG^{-2}\log T} (TK)^{-3/4}
(K^3\cdot K^{-1/2}T^{-1/2} + K)\cr& \ll_\e T^{\e}\max_{ K \le
TG^{-2}\log T} (K^{7/4}T^{-1/4} + T^{1/4}K^{1/4}) \cr& \ll_\e
T^{3/2+\e}G^{-7/2} + T^{1/2+\e} \ll_\e T^{1+\e}\cr} $$ for $G \ge
T^{1/7+\e}$. However, our initial condition was $G \ge
T^{1/5+\e}$, which is more restrictive. Fortunately, this is a
technical point that can be resolved by modifying Lemma 6
suitably. Namely, instead of $n^{3/2}T^{-1/2}$ in (4.4) we may put
$Cn^{5/2}T^{-3/2}$, which will be ``removed" in the fashion of
Lemma 6. Instead of the function $\a(u)$, another oscillating
function $\b(u)$ will appear, for which the analogue of (4.8) will
hold. Jutila obtained (4.5)-(4.6) by exploiting the fact that the
inversion made in (4.4) can be connected to the {\it Airy
integral} $$ {\roman {Ai}}(x) := {1\over\pi}\int_0^\infty
\cos({\txt{1\over3}}t^3 + tx)\d t \qquad(x \ge 0), $$ for which
there exist representations in terms of the classical Bessel
functions, thereby providing quickly asymptotic expansions
necessary in Lemma 6. In the new case there will be no Airy
integrals involved, but the necessary asymptotic expansion can be
obtained by the use of the saddle point method.

This ends the discussion of the case $m = 3$. The case $m=4$ will be
analogous, the non-trivial contribution will come from integer
quadruples $(n_1,n_2,n_3,n_4)$ such that $K < n_1,n_2,n_3,n_4
\le 2K $ $(T^{1/3} \le K \le TG^{-2}\log T)$ and
$$
|\sqrt{n_1} + \sqrt{n_2} - \sqrt{n_3} - \sqrt{n_4}\,|
\le T^{\e-1/2}.\leqno(4.12)
$$
Instead of Lemma 4 we use Lemma 5 (with $k=2$), to obtain
a contribution which is
$$
\eqalign{&
\ll_\e T^{1+\e}\max_{K\le TG^{-2}\log T}
(TK)^{-1}(K^4\cdot K^{-1/2}T^{-1/2} + K^2)\cr&
\ll_\e T^\e\max_{K\le TG^{-2}\log T}(T^{-1/2}K^{5/2} + K)\cr&
\ll_\e T^{2+\e}G^{-5} + T^{1+\e} \ll_\e T^{1+\e}\cr}
$$
for $G \ge T^{1/5+\e}$, as asserted. In this case direct application
of Lemma 6 suffices. Values of $m$ satisfying $m>4$ in (4.1) could
be handled in a similar fashion, provided that one can find analogues
of Lemma 4 and Lemma 5 which are strong enough.

\head
5. The moments of $J_2(t,G)$
\endhead
We shall prove now the analogue of Theorem 1 for $J_2(t,G)$. This is a more
difficult problem, and the ranges for $G$ for which the analogue
of (4.1) will hold will be poorer. The result is

\medskip
THEOREM 2. {\it We have}
$$
\int_T^{2T}J_2^m(t,G)\d t \ll_\e T^{1+\e}\leqno(5.1)
$$
{\it for $T^{1/2+\e} \le G \le T$ if $m = 1,2$; for
$T^{4/7+\e} \le G \le T$ if
$m=3$, and for $T^{3/5+\e} \le G \le T$ if $m=4$.}

\medskip
{\bf Proof of Theorem 1.} Our starting point is Lemma 3. We remark
that by using the estimate (see Y. Motohashi [22, Section 3.4])
$$
\sum_{\k_j\le K}\a_j\H \;\ll_\e\; K^2\log^3K\leqno(5.2)
$$
we see that for $G \ge T^{2/3}\log^CT$ the right-hand side of (3.4)
is $O(1)$, hence we may suppose that $T^{1/2+\e} \le G \le T^{2/3}\log^CT$.
Observe also that we may truncate the series in (3.4) at $TG^{-1}\log T$
with a negligible error. Besides (5.2) we need one more ingredient
from spectral theory, namely the author's bound [14]
$$
\sum_{K-1\le \k_j\le K+1}\a_j\H \ll_\e K^{1+\e}.\leqno(5.3)
$$
We give now the proof of Theorem when $m=2$ (the case $m=1$ easily
follows from this and the Cauchy-Schwarz inequality). We use (4.7)
with $J_2$ replacing $J_1$. Then
$$
\eqalign{&
\int_T^{2T}J_2^2(t,G)\d t \ll_\e T^{1+\e} + T^\e\max_{K \le TG^{-1}\log T}
\sum_{K<\k_j,\k_\ell\le2K}\a_j\a_l\H H_\ell^3(\hf)\times
\cr&
\times(\k_j\k_\ell)^{-1/2}\left(\k_j\over4{\roman e}\right)^{i\k_j}
\left(\k_\ell\over4{\roman e}\right)^{-i\k_\ell}
\int_T^{2T}\f(t)t^{i\k_j-i\k_\ell}\exp\left(-{\txt{1\over4}}G^2t^{-2}
(\k_j^2+\k_\ell^2)\right)\d t.\cr&
\cr}
$$
Repeated integrations by parts show that the contribution of
$\k_j,\,\k_\ell$ for which $|\k_j - \k_\ell| \ge T^\e$ is negligible. The
contribution of $|\k_j - \k_\ell| < T^\e$ is estimated by (5.3)
(splitting the summation over $\k_\ell$ in subsums of length $\le 2$)
and (5.2), while the integral over $t$ is estimated
trivially. The contribution will be
$$
\eqalign{&
\ll_\e T^{1+\e}\max_{K\le TG^{-1}\log T}T^{-1}\sum_{K<\k_j\le2K}\a_j\H
\k_j^{-1/2}\sum_{|\k_j-\k_\ell|<T^\e}\a_\l H_\ell^3(\hf)\k_\ell^{-1/2}\cr&
\ll_\e T^{\e}\max_{K\le TG^{-1}\log T}
\sum_{K<\k_j\le2K}\a_j\H K^{-1/2}K^{1/2}
\cr&
\ll_\e T^{2+\e}G^{-2} \ll_\e T^{1+\e}\cr}
$$
for $G \ge T^{1/2+\e}$, as asserted. We remark that the technique of this
proof can be used, following the arguments in [9, Chapter 5], to yield a
quick proof of the important bound
$$
\int_0^T E_2^2(t)\d t \ll_\e T^{2+\e},
$$
which is only slightly weaker than the second bound in (3.3).

The cases $m=3$ and $m=4$ are dealt with analogously. For the former,
after we raise the sum in (3.4) to the cube, it is seen that the
non-negligible contribution comes from the triplets
$(\k_j,\k_m,\k_\ell)$ for which $|\k_j + \k_m - \k_\ell | < T^\e$.
For the summation over one of the variables, say $\k_\ell$, we use
(5.3), and for the summation over $\k_j,\,\k_m$ we use (5.2). Similarly,
in the case of the fourth power the non-negligible contribution will
come from quadruples $(\k_j,\k_m,\k_\ell,\k_n)$ for which
$$
|\k_j + \k_m - \k_\ell- \k_n | < T^\e.\leqno(5.4)
$$
In this way the assertions of Theorem 2 concerning the cases $m = 3,4$
are obtained; the details are omitted for the sake of brevity.  Note
that $\lambda_j = \k_j^2 + {1\over4}$, so that (5.4) can be rewritten as
$$
|\sqrt{\lambda_j} + \sqrt{\lambda_m}- \sqrt{\lambda_\ell}-
\sqrt{\lambda_n}\,| < T^\e,\leqno(5.5)
$$
which is somewhat analogous to (4.12). Lemma 5 provides a good bound
for the number of integer quadruples satisfying (4.12), but the
condition (5.5) is much more difficult to deal with, since little is
known about arithmetic properties of the spectral values $\lambda_j$.

\head
4. Bounds for moments of $|\zt|$
\endhead
We shall show now how the results on power moments of $|\zt|$ follow
from mean square results on short intervals. In particular, a new
result will be derived, which connects power moments of $|\zt|$
with upper bounds furnished by Theorem 1 and Theorem 2.

\medskip
To begin with, suppose that $\{t_r\}_{r=1}^R$ are points lying in $[T,\,2T]$
such that $t_{r+1}-t_r\ge1\;(r=1,\ldots,R-1)$ and $|\z(\hf+it_r)| \ge
V \ge T^\e$ for $r =1,\ldots,R$. From (1.5) we have
$$
\eqalign{
RV^{2k} &\ll \log T\sum_{r=1}^R\int_{t_r-1/3}^{t_r+1/3}|\zt|^{2k}\d t\cr&
\ll \log T\sum_{s=1}^S\int_{\tau_s-G}^{\tau_s+G}|\zt|^{2k}\d t\cr&
\ll G\log T\sum_{s=1}^S J_k(\tau_s,G),\cr}\leqno(6.1)
$$
where we have grouped integrals over disjoint intervals
$[t_r-1/3,\,t_r+1/3]$ into integrals over disjoint intervals
$[\tau_s - G,\,\tau_s + G]$ with $s = 1,\ldots,S \;(\le R),\,G\ge T^\e$.

Suppose now that $k = 1$ in (6.1). Then we use Lemma 2, noting
that there are no absolute value signs on the right-hand side of
(2.9), which can be truncated at $TG^{-2}\log T$ and where, as
before, we may assume that $G \le T^{1/3}$. Exchanging the order
of summation, it follows from (6.1) that $$\eqalign{& RV^{2} \ll
RG\log^2T \;+\cr& + \sqrt{2}G\log T\sum_{n\le TG^{-2}\log
T}(-1)^nd(n)n^{-1/2}\times\cr&\times \sum_{s=1}^S\left(
\left({t_s\over2\pi n}+{1\over4}\right)-{1\over2}\right)^{-1/2}
\exp(-G^2\cdots)\sin f(t_s,n)\cr& \ll RG\log^2T +
G\Bigl(\sum_{n\le TG^{-2}\log T}d^2(n)n^{-1}\bigr)^{1/2}
\Bigl(\sum_{n\le TG^{-2}\log T}
\Bigl|\sum_{s=1}^S\cdots\Bigr|^2\Bigr)^{1/2},\cr}\leqno(6.2) $$ by
the Cauchy-Schwarz inequality. We further have $$ \sum_{n\le
TG^{-2}\log T} \Bigl|\sum_{s=1}^S\cdots\Bigr|^2 = \sum_{s_1,s_2\le
S}S_0, $$ where $$\eqalign{ S_0 &:= \sum_{n\le TG^{-2}\log T}
\left(\left({t_{s_1}\over2\pi n}+{1\over4}\right)^{1/2} -
{1\over2}\right)^{-1/2}\cr& \times\exp(-G^2({\roman
{arsinh}}\sqrt{\pi n/(2t_{s_1})})^2) \left(\left({t_{s_2}\over2\pi
n}+{1\over4}\right)^{1/2} - {1\over2}\right)^{-1/2}\cr&\times
\exp(-G^2({\roman {arsinh}}\sqrt{\pi n/(2t_{s_2})})^2)
\exp(if(t_{s_1},n) - if(t_{s_2},n)).\cr} $$ Removing by partial
summation monotonic coefficients from $S_0$, we are led to the
crucial exponential sum $$ S_1 : = \sum_{n\le M}
\exp(if(t_{s_1},n) - if(t_{s_2},n)) \qquad(M \le TG^{-2}\log
T).\leqno(6.3) $$ The quality of the estimation of $S_1$ is
limited by the scope of the present-day exponential sum techniques
(see e.g., M.N. Huxley [4]). The terms $s_1 = s_2$ in (6.3) will
eventually give rise to $R \ll_\e T^{1+\e}V^{-6}$, namely to a
weak form of the sixth moment (3.6), but it does not seem likely
that (3.6) can be reached (unconditionally) in this fashion.
Observing that
$$
 {\partial f(x,k)\over\partial x} \,=\,
2\,\text{arsinh}\,\sqrt{\pi k \over2x} \;\sim\; \sqrt{2\pi k\over
x}\quad(x \asymp T,\,k \le TG^{-2}\log T),\leqno(6.4) $$ setting
$$ f(u) := f(t_{s_1},u) - f(t_{s_2},u), \quad F :=
|t_{s_1}-t_{s_2}| (KT)^{-1/2}, $$ we have (see [8, Chapters 1-2]
for the relevant exponent pair technique) that $S_1$ may be split
into $O(\log T)$ subsums of the type $$\eqalign{& \sum_{K<k\le
K'\le2K}\exp(if(k)) \ll F^\kappa K^\lambda + F^{-1}\cr& \ll
J^{\kappa}T^{-\kappa/2}K^{\lambda-\kappa/2} + (KT)^{1/2}
|t_{s_1}-t_{s_2}|^{-1}\quad(K \le TG^{-2}\log T),\cr} $$ provided
that $|t_{s_1}-t_{s_2}| \le J (\ll T)$, and $(\kappa,\lambda)$ is
a (one-dimensional) exponent pair. Choosing $(\kappa,\lambda) =
(\hf,\hf)$, $J = T^{-\e}G^3$ we obtain (2.7) (with $T^\e$ in place
of $\log^{17}T$). Namely with $J = T^{-\e}G^3$ the number of
points $R = R_0$ to be estimated satisfies $R_0 \ll_\e
T^{1+\e}G^{-3}$, hence dividing $[T/2,\,5T/2]$ into subintervals
of length not exceeding $J$ one obtains $$ R \;\ll\; R_0(1 + T/J)
\;\ll_\e\; T^{2+\e}G^{-6} \;\ll_\e\; T^{2+\e}V^{-12}, $$ which
easily yields (2.7) (with $T^\e$ in place of $\log^{17}T$). This
analysis was carried in detail in [8, Chapter 8], where the
possibilities of choosing other exponent pairs besides $(\kappa,
\lambda) = (\hf,\hf)$  were discussed.

Another type of a similar estimate was obtained by the author
in [15] (for the analysis
of sums of moments over well-spaced points $\{t_r\} \in [T,\,2T]$ the
reader is referred to [10]). This result will be stated here as

\medskip
THEOREM 3. {\it Let $T \le t_1 < \ldots < t_R \le 2T$ be points such that
$|\z(\hf + it_r)| \ge VT^{-\e}$ with
$t_{r+1} - t_r \ge V \ge T^{{1\over10}+\e}$ for $r = 1,\ldots, R-1$.
Then, for any fixed integer $M\ge 1$,}
$$\eqalign{&
R \ll_\e T^{\e-M/2}V^{-2}\max_{K\le T^{1+\e}V^{-4}}\times\cr&\times
\int_{T/2}^{5T/2}\f(t)\Big|\sum_{K\le k\le K'\le2K}(-1)^kd(k)k^{-1/4}
\exp(2i\sqrt{2\pi kt}+cik^{3/2}t^{-1/2})\Big|^{2M}\d t,\cr}\leqno(6.4)
$$
{\it where $c = \sqrt{2\pi^3}/6$ and $\f(t)$ is a non-negative, smooth
function supported in $[T/2,5T/2]$ such that $\f(t) = 1$ for $T\le t \le2T$}.

\medskip
The case $M = 1$ quickly leads to a weakened form of the fourth moment
estimate, namely $\int_0^T|\zt|^4\d t \ll_\e T^{1+\e}$. The case $M = 2$
of (6.4),
by the use of Lemma 5 and Lemma 6, will lead again to a weakened form
of the twelfth moment bound (2.7) (with $T^{2+\e}$).

\medskip
Finally we present a new result, which connects bounds for moments of
$|\zt|$ to bounds of moments of $J_k(t,G)$. This is

\medskip
THEOREM 4. {\it Suppose that}
$$
\int_T^{2T}J_k^m(t,G)\d t \ll_\e T^{1+\e}\leqno(6.5)
$$
{\it holds for some fixed $k,m  \in\NN$ and $G \ge T^{\a_{k,m}+\e},\,
 0 \le \a_{k,m} < 1$. Then}
$$
\int_0^T|\zt|^{2km}\d t \ll_\e T^{1+(m-1)\a_{k,m}+\e}.\leqno(6.6)
$$

\medskip
{\bf Proof of Theorem 4.} We note first that ($\mu(\cdot)$ denotes
measure) the bound
$$ \mu\Bigl(t\in[T,\,2T] \,:\,J_k(t,G)\ge
U\Bigr) \ll_\e T^{1+\e}U^{-m} \leqno(6.7) $$ follows from (6.5).
We use (6.1), dividing the sum over $s$ into $O(\log T)$ subsums
where $U < J_k(\tau_s,G) \le 2U$. Then, for $U_0\;(\gg1)$ to be
determined later, $$ \eqalign{ &\sum_{s=1}^S J_k(\tau_s,G) \ll
SU_0 + \log T\max_{U\ge U_0} \sum_{s,U < J_k(\tau_s,G) \le
2U}J_k(\tau_s,G)\cr& \ll SU_0 + U\log T\max_{U\ge U_0} \sum_{s,U <
J_k(\tau_s,G) \le 2U}1\cr& \ll_\e SU_0 + G^{-1}\log T\max_{U\ge
U_0} T^{1+\e}U^{1-m}\cr& \ll_\e SU_0 +
T^{1+\e}U_0^{1-m}G^{-1},\cr} $$ since $m\ge1$ and $$
\sum_{s,J_k(\tau_s,G)>U}1 \ll_\e T^{1+\e}U^{1-m}G^{-1}.\leqno(6.8)
$$ Namely if $J_k(\tau_s,G) > U$, then $J_k(t,2G) \ge U$ for
$t\in[\tau_s-\hf G,\,\tau_s + \hf G]$, and (6.8) follows from
(6.7). The choice $$ U_0 \;=\; {\left({T\over
SG}\right)}^{1/m}\quad(\,\gg 1\,) $$ yields $$ \sum_{s=1}^S
J_k(\tau_s,G) \ll_\e T^{1/m+\e}S^{1-1/m}G^{-1/m}.\leqno(6.9) $$
Inserting (6.9) in (6.1) we obtain $$ R \ll_\e
T^{1+\e}G^{m-1}V^{-2km},\leqno(6.10) $$ and (6.6) easily follows
from (6.10), on taking $G = T^{\a_{k,m}+\e}$.

\medskip
This completes the proof of Theorem 4. The values $\a_{1,2} = 0$
(Theorem 1) and $\a_{2,2} = \hf$ (Theorem 2) yield, respectively,
$$ \int_0^T|\zt|^4\d t \ll_\e T^{1+\e},\quad \int_0^T|\zt|^8\d t
\ll_\e T^{3/2+\e}.\leqno(6.11) $$ The bounds in (6.11) are, of
course, well-known, but they are (up to the factor $T^\e$) the
sharpest known ones, and the bound for the fourth moment is
essentially of the correct order of magnitude. Other values of
$\a_{k,m}\;(k=1,2)$, furnished by Theorem 1 and Theorem 2, do not
yield any new bounds, as can be readily checked. However, it seems
that this approach is of interest, especially in view of recent
results on the distribution of sums and differences of square
roots of integers (cf. Lemma 4 and Lemma 5).

\vfill
\eject
\topglue2cm
\Refs
\bigskip

\item{[1]} F.V. Atkinson, The mean value of the Riemann zeta-function,
Acta Math. {\bf81}(1949), 353-376.

\item{[2]} J.B. Conrey, D.W. Farmer, J.P. Keating, M.O. Rubinstein
and N.C. Snaith, Integral moments of $L$-functions, preprint, 58pp,
arXiv:math.NT/0206018,

{\tt http://front.math.ucdavis.edu/mat.NT/0206018}.

\item{[3]} J.L. Hafner and A. Ivi\'c, On the mean square of the
Riemann zeta-function on the critical line.  J. Number Theory,
{\bf 32}(1989), 151--191.

\item{[4]} D.R. Heath-Brown, The twelfth power moment of the Riemann
zeta-function, Quart. J. Math. (Oxford) {\bf29}(1978), 443-462,

\item{[5]}  D.R. Heath-Brown, The fourth moment of the Riemann
zeta-function,
 Proc. London Math. Soc.  (3){\bf 38}(1979), 385-422.

\item{[6]} M.N. Huxley, Area, Lattice Points and Exponential
Sums, Oxford Science Publications, Clarendon Press,
Oxford, 1996

\item{[7]}  A.E. Ingham, Mean-value theorems in the theory of the Riemann
zeta-function,  Proc. London Math. Soc.  (2){\bf27}(1926), 273-300.

\item{[8]} A. Ivi\'c, The Riemann zeta-function, John Wiley \&
Sons, New York, 1985.

\item{[9]} A. Ivi\'c, The mean values of the Riemann zeta-function,
LNs {\bf 82}, Tata Inst. of Fundamental Research, Bombay (distr. by
Springer Verlag, Berlin etc.), 1991.

\item{[10]} A. Ivi\'c, Power moments of the Riemann zeta-function
over short intervals, Arch. Mat. {\bf62}\ (1994),\ 418-424.

\item{[11]} A. Ivi\'c,  On the fourth moment of the Riemann
zeta-function, Publs. Inst. Math. (Belgrade) {\bf 57(71)}(1995), 101-110.

\item{[12]} A. Ivi\'c,
 On the ternary additive divisor problem and the sixth moment of the
zeta-function, ``Sieve Methods, Exponential Sums, and their Applications in
Number Theory" (eds. G.R.H. Greaves, G. Harman, M.N. Huxley), Cambridge
University Press, Cambridge, 1996, 205-243.

\item{[13]} A. Ivi\'c, On some results concerning the Riemann Hypothesis, in
``Analytic Number Theory", LMS LNS 247, Cambridge University Press,
Cambridge, 1997, pp. 139-167.

\item{[14]} A. Ivi\'c,  On sums of Hecke series in short intervals,
J. de Th\'eorie des Nombres Bordeaux {\bf13}(2001), 453-468.

\item{[15]} A. Ivi\'c, Sums of squares of $|\zt|$ over short
intervals, Max-Planck-Institut f\"ur Mathematik, Preprint Series
2002({\bf52}), 12 pp. arXiv:math.NT/0311516,

{\tt http://front.math.ucdavis.edu/mat.NT/0311516}.

\item{[16]} A. Ivi\'c and Y. Motohashi, The mean square of the
error term for the fourth moment of the zeta-function,  Proc. London Math.
Soc. (3){\bf66}(1994), 309-329.

\item {[17]}  A. Ivi\'c and Y. Motohashi,  The fourth moment of the
Riemann zeta-function,  J. Number Theory  {\bf51}(1995), 16-45.

\item{[18]} A. Ivi\'c and P. Sargos, On the higher moments of the
error term in the divisor problem, to appear.

\item{[19]} M. Jutila, Riemann's zeta-function and the divisor problem,
Arkiv Mat. {\bf21}(1983), 75-96 and II, ibid. {\bf31}(1993), 61-70.

\item{[20]} K. Matsumoto, Recent developments in the mean square theory
of the Riemann zeta and other zeta-functions, in ``Number Theory",
Birkh\"auser, Basel, 2000, 241-286.

\item {[21]}  Y. Motohashi, A relation between the Riemann zeta-function
and the hyperbolic Laplacian,  Ann. Sc. Norm. Sup. Pisa, Cl. Sci. IV
ser.  {\bf22}(1995), 299-313.

\item {[22]}  Y. Motohashi,  Spectral theory of the Riemann
zeta-function,  Cambridge University Press, Cambridge, 1997.

\item{[23]} K. Ramachandra, A simple proof of the mean fourth
power estimate for $L({1\over2}+it,\chi)$,  Ann. Sc. Norm. Sup. Pisa,
Cl. Sci. {\bf1}(1974), 81-97.

\item{[24]} K. Ramachandra, On the mean-value and omega-theorems
for the Riemann zeta-function, Tata Institute of Fundamental Research,
Bombay, distr. by Springer Verlag, 1995.

\item{[25]} O. Robert and P. Sargos, Three-dimensional
exponential sums with monomials, J. reine angew. Math. (in print).

\item{[26]} E.C. Titchmarsh, The theory of the Riemann zeta-function
(2nd ed.),  University Press, Oxford, 1986.

\vskip2cm
\endRefs
{\sevenrm Aleksandar Ivi\'c \par Katedra Matematike RGF-a\par
Universiteta u Beogradu\par \DJ u\v sina 7, 11000 Beograd\par
Serbia and Montenegro}\par {\sevenbf e-mail: aivic\@matf.bg.ac.yu,
aivic\@rgf.bg.ac.yu}

\bye